\documentclass[8pt]{article}%
\usepackage{amssymb}
\usepackage{amsmath}
\usepackage{amsfonts}
\usepackage{graphicx}%
\providecommand{\U}[1]{\protect \rule{.1in}{.1in}}
\newtheorem{theorem}{Theorem}

\newtheorem{proposition}[theorem]{Proposition}

\begin{document}

\title{On the Laplace Normal Vector Field \\of Skew Ruled Surfaces}
\author{Stylianos Stamatakis\textbf{\medskip}\\ \emph{Aristotle University of Thessaloniki}\\ \emph{Department of Mathematics}\\ \emph{GR-54124 Thessaloniki, Greece}\\ \emph{e-mail: stamata@math.auth.gr}}
\date{}
\maketitle

\begin{abstract}
We consider the Laplace normal vector field of relatively normalized ruled
surfaces with non-vanishing Gaussian curvature in the Euclidean space
$\mathbb{R}^{3}$. We determine all ruled surfaces and all relative normalizations for which
the Laplace normal image degenerates into a point or into a curve.
Moreover, we study the Laplace normal image of a
non-conoidal ruled surface whose relative normals lie on the asymptotic plane. \smallskip

 \smallskip


\end{abstract}

\section{Introduction}

Relatively normalized ruled surfaces with non-vanishing Gaussian curvature in
the Euclidean space $\mathbb{R}^{3}$ have been studied in the last years in
many points of view (see \cite{Heil}, \cite{Manhart1}, \cite{Manhart2},
\cite{Manhart3}; for more details and references see
\cite{Stamouetc}). This paper deals with the Laplace normal vector field of
skew ruled surfaces. We first show that for any relative normalization
$\bar{y}$ the Laplace normal vector of a ruled surface $\Phi$ along each
ruling lies on the corresponding asymptotic plane. Then, we determine all ruled
surfaces and all relative normalizations, for which the Laplace normal image
degenerates into a point or into a curve. We finish the paper by the study
of the Laplace normal image of a non conoidal ruled surface such relatively
normalized, that the relative normals $\bar{y}$ along each ruling lie on the
corresponding asymptotic plane.

\section{Preliminaries}

To set the stage for this work the classical notations of relative
differential geometry and of ruled surfaces theory are briefly presented. For
this purpose the books \cite{Hoschek} and \cite{Schirokow} are used as general
references. \medskip

\noindent In the Euclidean space $%
%TCIMACRO{\U{211d} }%
%BeginExpansion
\mathbb{R}
%EndExpansion
^{3}$ let $\Phi$ be a skew (non-developable)\ ruled $C^{r}$-surface, $r\geq2$.
We denote by $\bar{s}(u),$ $u\in I$ ($I\subset%
%TCIMACRO{\U{211d} }%
%BeginExpansion
\mathbb{R}
%EndExpansion
$ open interval) the position vector of the line of striction of $\Phi$\ and
by $\bar{e}(u)$ the unit vector pointing along the rulings. We choose the
parameter $u$ to be the arc length along the spherical curve $\bar{e}(u)$.
Then a parametrization of the ruled surface $\Phi$ over the region $U:=I\times%
%TCIMACRO{\U{211d} }%
%BeginExpansion
\mathbb{R}
%EndExpansion
$ of the $\left(  u,v\right)  $-plane is%
\begin{equation}
\bar{x}(u,v)=\bar{s}(u)+v\bar{e}(u)\quad \text{with}\quad \left \vert \bar
{e}\right \vert =|\bar{e}\,%
%TCIMACRO{\U{b4}}%
%BeginExpansion
\acute{}%
%EndExpansion
\,|=1,\text{\quad}\langle \bar{s}\,%
%TCIMACRO{\U{b4}}%
%BeginExpansion
\acute{}%
%EndExpansion
\,(u),\bar{e}\,%
%TCIMACRO{\U{b4}}%
%BeginExpansion
\acute{}%
%EndExpansion
\,(u)\rangle=0\text{ in }I, \label{1.20}%
\end{equation}
where the prime denotes differentiation with respect to $u$ and $\langle
~,~\rangle$\ the standard scalar product in $%
%TCIMACRO{\U{211d} }%
%BeginExpansion
\mathbb{R}
%EndExpansion
^{3}$. The \textit{Kruppa frame }of $\Phi,$ consisting of the vector $\bar
{e}(u)$, the central normal vector $\bar{n}:=\bar{e}$\thinspace$%
%TCIMACRO{\U{b4}}%
%BeginExpansion
\acute{}%
%EndExpansion
$ and the central tangent vector $\bar{z}:=\bar{e}\times \bar{n},$ satisfies
the relations \cite[p. 17]{Hoschek}%
\begin{equation}
\bar{e}\,%
%TCIMACRO{\U{b4}}%
%BeginExpansion
\acute{}%
%EndExpansion
=\bar{n},\quad \bar{n}\,%
%TCIMACRO{\U{b4}}%
%BeginExpansion
\acute{}%
%EndExpansion
=-\bar{e}+\kappa \bar{z},\quad \bar{z}\,%
%TCIMACRO{\U{b4}}%
%BeginExpansion
\acute{}%
%EndExpansion
=-\kappa \bar{n}, \label{1.21}%
\end{equation}
where
\[
\kappa=(\bar{e},\bar{e}\,%
%TCIMACRO{\U{b4}}%
%BeginExpansion
\acute{}%
%EndExpansion
\acute{} ,\bar{e}\, \acute{} \, \acute{} \,)
\]
denotes the \textit{conical curvature }of $\Phi.$ Consider the
\textit{parameter of distribution }
\[
\delta=(\bar{s}\,%
%TCIMACRO{\U{b4}}%
%BeginExpansion
\acute{}%
%EndExpansion
,\bar{e},\bar{e}\,%
%TCIMACRO{\U{b4}}%
%BeginExpansion
\acute{}%
%EndExpansion
\,)
\]
and the \textit{striction}
\[
\sigma:=\sphericalangle(\bar{e},\bar{s}\,%
%TCIMACRO{\U{b4}}%
%BeginExpansion
\acute{}%
%EndExpansion
\,)\quad(-\frac{\pi}{2}<\sigma \leq \frac{\pi}{2},\operatorname*{sign}%
\sigma=\operatorname*{sign}\delta)
\]
of $\Phi$. Then, the tangent vector $\bar{s}\,%
%TCIMACRO{\U{b4}}%
%BeginExpansion
\acute{}%
%EndExpansion
$ of the line of striction and the unit normal vector $\bar{\xi}$ of $\Phi$
are expressed by
\begin{equation}
\bar{s}\,%
%TCIMACRO{\U{b4}}%
%BeginExpansion
\acute{}%
%EndExpansion
=\delta \left(  \lambda \bar{e}+\bar{z}\right)  , \label{1.23}%
\end{equation}%
\begin{equation}
\bar{\xi}=\frac{\delta \bar{n}-v\bar{z}}{w}, \label{1.24}%
\end{equation}
where
\[
\lambda:=\cot \sigma
\]
and
\[
w^{2}:=v^{2}+\delta^{2}.
\]
When the \textit{fundamental invariants} $\kappa(u),\delta(u)$ and
$\lambda(u)$ are given, then there exists up to rigid motions of the space $%
%TCIMACRO{\U{211d} }%
%BeginExpansion
\mathbb{R}
%EndExpansion
^{3}$ a unique ruled surface $\Phi$, whose fundamental invariants are the
given. \smallskip

\noindent The components $h_{ij}$ of the second fundamental form in the coordinates
$u^{1}:=u$, $u^{2}:=v$ are the following%
\begin{equation}
h_{11}=-\frac{\kappa v^{2}{\small +}\delta \,%
%TCIMACRO{\U{b4}}%
%BeginExpansion
\acute{}%
%EndExpansion
\,v{\small +}\delta^{2}\left(  \kappa-\lambda \right)  }{w},\quad h_{12}%
=\delta,\quad h_{22}=0. \label{1.25}%
\end{equation}
The Gaussian curvature $K$ of $\Phi$\ is given by%
\begin{equation}
K=\frac{-\delta^{2}}{w^{4}}. \label{1.27}%
\end{equation}

A $C^{s}$-mapping $\bar{y}:U\longrightarrow%
%TCIMACRO{\U{211d} }%
%BeginExpansion
\mathbb{R}
%EndExpansion
^{3},$ $r>s,$ is called a $C^{s}$-\textit{relative normalization} of $\Phi$
if
\begin{equation}
\operatorname*{rank}\left \{  \bar{x}_{/1},\bar{x}_{/2},\bar{y}\right \}
=3,\quad \operatorname*{rank}\left \{  \bar{x}_{/1},\bar{x}_{/2},\bar{y}%
_{/1},\bar{y}_{/2}\right \}  =2\quad \text{ }\forall \text{ }\alpha \in
U.~\footnote{Partial derivatives of a function $f$ are denoted by
$f_{/i}:=\frac{\partial f}{\partial u^{i}},~f_{/ij}:=\frac{\partial^{2}%
f}{\partial u^{i}\partial u^{j}}$ etc. $\medskip$} \label{1.1a}%
\end{equation}
The\  \textit{covector }$\bar{X}$ of the tangent plane is defined by%
\begin{equation}
\langle \bar{X},\bar{x}_{/1}\rangle=\langle \bar{X},\bar{x}_{/2}\rangle
=0\quad \text{and}\quad \langle \bar{X},\bar{y}\rangle=1. \label{1.2}%
\end{equation}
The \textit{relative metric} $G$ on $U$ is introduced by
\begin{equation}
G_{ij}=\langle \bar{X},\bar{x}_{/ij}\rangle. \label{1.3}%
\end{equation}
The \textit{support function }of the relative normalization $\bar{y}$ is
defined, according to \cite{Manhart3}, by%
\[
q:=\langle \bar{\xi},\bar{y}\rangle:U\longrightarrow%
%TCIMACRO{\U{211d} }%
%BeginExpansion
\mathbb{R}
%EndExpansion
,\quad q\in C^{s}\left(  U\right)  ,
\]
where $\bar{\xi}:U\longrightarrow%
%TCIMACRO{\U{211d} }%
%BeginExpansion
\mathbb{R}
%EndExpansion
^{3}$ is the \textit{Euclidean normalization} of $\Phi$. By virtue of
(\ref{1.1a}), $q$ never vanishes on $U;$ moreover, because of (\ref{1.2}), it
turns out
\begin{equation}
\bar{X}=q^{-1}\bar{\xi}. \label{1.4}%
\end{equation}
On account of (\ref{1.3}) and (\ref{1.4}) we obtain%
\begin{equation}
G_{ij}=q^{-1}h_{ij}. \label{1.5}%
\end{equation}
We mention that given a support function $q,$ the relative normalization
$\bar{y}$ is uniquely determined and possesses the following parametrization
\cite[p. 197]{Manhart3}
\begin{equation}
\bar{y}=-h^{ij}\,q_{/i}\, \bar{x}_{/j}+q\, \bar{\xi}, \label{1.6}%
\end{equation}
where $h^{ij}$ are the components of the inverse tensor of $h_{ij}$.
\\

In \cite{Heil}, H. Heil has introduced the $C^{1}$-mapping $\bar{L}:U\longrightarrow%
%TCIMACRO{\U{211d} }%
%BeginExpansion
\mathbb{R}
%EndExpansion
^{3}$, defined by%
\begin{equation}
\bar{L}(u,v)=\frac{\triangle \bar{x}(u,v)}{2}, \label{1.6c}%
\end{equation}
where $\triangle$ is the laplacian with respect to the relative metric $G,$ as
\textit{the Laplace normal vector }of $\Phi$. When we move the
Laplace normal vectors to a fixed point, the endpoints of them describe the
\textit{Laplace normal image of }$\Phi.$

%-----------------------------------------------------------------------------------------------------------------

\section{The Laplace normal vector field of a ruled surface}

We consider a relatively normalized skew ruled $C^{2}$-surface $\Phi \subset%
%TCIMACRO{\U{211d} }%
%BeginExpansion
\mathbb{R}
%EndExpansion
^{3}$. By virtue of (\ref{1.20})-(\ref{1.25}) and (\ref{1.6}), the relative
normalization $\bar{y}$ can be written as follows%
\begin{equation}
\bar{y}=-w\frac{\delta q_{/1}+q_{/2}(\kappa w^{2}+\delta \,%
%TCIMACRO{\U{b4}}%
%BeginExpansion
\acute{}%
%EndExpansion
\,v)}{\delta^{2}}\bar{e}+\frac{\delta^{2}q-w^{2}vq_{/2}}{\delta w}\bar
{n}-\frac{vq+w^{2}q_{/2}}{w}\bar{z}. \label{1.6b}%
\end{equation}
We denote by
\begin{equation}
q_{_{AFF}}:=|K|^{1/4} \label{1.6a}%
\end{equation}
the support function of the \textit{equiaffine normalization }$\bar{y}%
_{_{AFF}}$ \textit{of }$\Phi.$ On account of (\ref{1.27}), (\ref{1.6b}) and
(\ref{1.6a}) we get%
\begin{equation}
\bar{y}_{_{AFF}}=\frac{\varepsilon}{|\delta|^{1/2}}\left(  \frac{2\kappa
v+\delta \,%
%TCIMACRO{\U{b4}}%
%BeginExpansion
\acute{}%
%EndExpansion
\,}{2\delta}\bar{e}+\bar{n}\right)  , \label{1.7}%
\end{equation}
where $\varepsilon=\operatorname*{sign}\delta.$ Using (\ref{1.20}%
)-(\ref{1.25}), (\ref{1.5}) and (\ref{1.6c}) we find%

\begin{equation}
\bar{L}=\frac{wq}{\delta}\left(  \frac{2\kappa v+\delta \,%
%TCIMACRO{\U{b4}}%
%BeginExpansion
\acute{}%
%EndExpansion
}{2\delta}\bar{e}+\bar{n}\right)  , \label{2.1}%
\end{equation}
i.e.
\[
\bar{L}=\frac{q}{q_{_{AFF}}}\bar{y}_{_{AFF}}%
\]
which shows, that \textit{for any normalization }$\bar{y}$ \textit{the Laplace
normal $\bar{L}$ of }$\Phi$ \textit{along each ruling lies on the corresponding
asymptotic plane in the direction of the equiaffine normalization $\bar{y}_{_{AFF}}$ and
therefore it is independent of the relative normalization} $\bar{y}$.\medskip

%------------------------------------------------------------------------------------------

\noindent \textbf{Remark}. Obviously, \textit{two ruled surfaces parametrized by
}(\ref{1.20}) \textit{with parallel rulings, common parameter of distribution
and common support function have the same Laplace normal vector field}%
.\medskip

We shall firstly determine all ruled surfaces, whose Laplace normals are
constant along each ruling. We put
\begin{equation}
L_{1}=\frac{wq\left(  2\kappa v+\delta \,%
%TCIMACRO{\U{b4}}%
%BeginExpansion
\acute{}%
%EndExpansion
\, \right)  }{2\delta^{2}},\quad L_{2}=\frac{wq}{\delta}, \label{2.5}%
\end{equation}
and on account of (\ref{1.21}) we find%
\begin{equation}
\bar{L}_{/1}=\left(  L_{1/1}-L_{2}\right)  \bar{e}+\left(  L_{1}%
+L_{2/1}\right)  \bar{n}+\kappa L_{2}\bar{z},\quad \bar{L}_{/2}=L_{1/2}\bar
{e}+L_{2/2}\bar{n}. \label{2.7}%
\end{equation}
Thus we have $\bar{L}_{/2}=\bar0$ if and only if $\kappa$ vanishes and the
support function is of the form%
\begin{equation}
q=\frac{f(u)}{w}, \label{2.7a}%
\end{equation}
where $f(u)$ is an arbitrary nonvanishing $C^{2}$-function. So we have the

%-----------------------------------------------------------------------------------------------

\begin{proposition}
The Laplace normals of a relatively normalized ruled $C^{3}$-surface $\Phi$,
which is free of torsal rulings, are constant along each ruling if and only if
$\Phi$ is conoidal and the support function is of the form (\ref{2.7a}).
\end{proposition}

In this case the Laplace normal becomes%
\begin{equation}
\bar{L}(u)=\frac{f}{\delta}\left(  \frac{\delta \,%
%TCIMACRO{\U{b4}}%
%BeginExpansion
\acute{}%
%EndExpansion
}{2\delta}\bar{e}+\bar{n}\right)  , \label{2.8}%
\end{equation}
i.e., \textit{the Laplace normal image of }$\Phi$\textit{ is degenerated into
a point or into a curve}. In the following we investigate all ruled surfaces,
which possess the previous property.\medskip

\noindent \textit{Case I}. The Laplace normal image degenerates into a point if and only
if $\operatorname*{rank}(\bar{L}_{/1},\bar{L}_{/2})=0,$ or, on account of
(\ref{2.7}), equivalently,%
\begin{equation}
L_{1/1}-L_{2}=L_{1}+L_{2/1}=\kappa~L_{2}=L_{1/2}=L_{2/2}=0.\label{2.10}%
\end{equation}
Thus $\kappa=0$ and the support function is of the form (\ref{2.7a}). In view
of (\ref{2.5})\ the remaining relations of (\ref{2.10}) become%
\begin{equation}
\left(  \frac{\delta \, \acute{}\,f}{\delta^{2}}\right)
%TCIMACRO{\QATOP{\acute{}}{{}}}%
%BeginExpansion
\genfrac{}{}{0pt}{}{\acute{}}{{}}%
%EndExpansion
-2\frac{f}{\delta}=0,\quad \frac{\delta \,%
%TCIMACRO{\U{b4}}%
%BeginExpansion
\acute{}%
%EndExpansion
\,f}{\delta^{2}}+2\left(  \frac{f}{\delta}\right)
%TCIMACRO{\QATOP{\acute{}}{{}}}%
%BeginExpansion
\genfrac{}{}{0pt}{}{\acute{}}{{}}%
%EndExpansion
=0.\label{2.13}%
\end{equation}
whence

\[
\left(  \frac{f}{\delta}\right)
\genfrac{}{}{0pt}{}{\acute{}~\acute{}}{{}}%
%EndExpansion
+\frac{f}{\delta}=0.
\]
Consequently
\[
\frac{f}{\delta}=c_{1}\cos u+c_{2}\sin u,
\]
where $c_{1},c_{2}=const.$ Then, the second relation of (\ref{2.13}) implies%
\[
\frac{\delta \,%
%TCIMACRO{\U{b4}}%
%BeginExpansion
\acute{}%
%EndExpansion
}{2\delta}=\frac{c_{1}\sin u-c_{2}\cos u}{c_{1}\cos u+c_{2}\sin u},
\]
from which we find
\[
\delta=c_{3}\left(  c_{1}\cos u+c_{2}\sin u\right)  ^{-2},
\]
where $c_{3}=const.\neq0.$ Thus
\[
f=\pm \left(  c_{3}\delta \right)  ^{1/2}.
\]
On account of (\ref{1.27}), (\ref{1.6a}) and (\ref{2.7a}) we have
\[
q=\pm|c_{3}|^{1/2}\,q_{_{AFF}}.
\]
Finally from (\ref{1.6b}) and (\ref{2.8})\ we find%
\[
\bar{y}=\bar{L}=\bar{a},
\]
where $\bar{a}$ is the constant vector%
\[
\bar{a}=\left(  c_{1}\sin u-c_{2}\cos u\right)  \bar{e}+\left(  c_{1}\cos
u+c_{2}\sin u\right)  \bar{n}.
\]
Therefore $\Phi$ is an improper relative sphere, see \cite{Blaschke},
\cite{Schirokow}. So we have

%-----------------------------------------------------------------------------------------------------------

\begin{proposition}
The Laplace normal image of a relatively normalized ruled $C^{3}$-surface
$\Phi$, which is free of torsal rulings, degenerates into a point if and only
if $\Phi$ is an improper conoidal relative sphere, the parameter of
distribution is given by $\delta=c_{3}\left(  c_{1}\cos u+c_{2}\sin u\right)
^{-2},$ where $c_{1},c_{2},c_{3}=const.,~c_{3}\neq0,$ and the support function
is given by $q=\pm|c_{3}|^{1/2}\,q_{_{AFF}}$.
\end{proposition}

\noindent \textit{Case II}.\textbf{ }The Laplace normal image degenerates into a curve
$\Gamma$ if and only if $\operatorname*{rank}(\bar{L}_{/1},\bar{L}_{/2})=1,$
or, on account of (\ref{2.7}) equivalently,%
\begin{equation}
\frac{L_{1/1}-L_{2}}{L_{1/2}}=\frac{L_{1}+L_{2/1}}{L_{2/2}}=\frac
{\kappa \,L_{2}}{0}. \label{2.14}%
\end{equation}
Obviously $\kappa=0.\medskip$

\noindent \textit{Subcase IIa. }Let $L_{2/2}=0.$ The support function has therefore the
form (\ref{2.7a}) and the curve $\Gamma$ has the parametrization (\ref{2.8}).
Moreover it is easy to confirm, that the curve $\Gamma$ is planar.\smallskip

\noindent \textit{Subcase IIa. }For $L_{2/2}\neq0$ it turns out from (\ref{2.14})%
\begin{equation}
2\delta \delta \,%
%TCIMACRO{\U{b4}}%
%BeginExpansion
\acute{}%
%EndExpansion
\,%
%TCIMACRO{\U{b4}}%
%BeginExpansion
\acute{}%
%EndExpansion
\,-3\delta \,%
%TCIMACRO{\U{b4}}%
%BeginExpansion
\acute{}%
%EndExpansion
\,^{2}-4\delta^{2}=0, \label{2.16}%
\end{equation}
from which we obtain%
\begin{equation}
\delta=c_{2}\cos^{-2}(u+c_{1}),\quad c_{1},c_{2}=const. \label{2.19a}%
\end{equation}
Hence, a parametrization of the curve $\Gamma$ is%
\begin{equation}
\bar{L}=\frac{wq}{c_{2}}\cos(u+c_{1})\  \bar{\alpha}, \label{2.20}%
\end{equation}
where $\bar{\alpha}$ is the constant unit vector%
\[
\bar{\alpha}=\sin(u+c_{1})\  \bar{e}+\cos(u+c_{1})\  \bar{n}.
\]
Obviously, $\Gamma$ is a straight line. Finally, from (\ref{1.7}) it follows%
\[
\bar{y}_{_{AFF}}=\pm \left \vert c_{2}\right \vert ^{-1/2}\bar{a},
\]
so that the affine normal image of $\Phi$ is an improper affine
sphere.\medskip

\noindent Thus we have proved

%-----------------------------------------------------------------------------------------------------------

\begin{proposition}
The Laplace normal image of a relatively normalized ruled $C^{3}$-surface
$\Phi$, which is free of torsal rulings, degenerates into a curve $\Gamma$ if
and only if $\Phi$ is conoidal and either the support function has the form
$q=f(u)w^{-1},$ or the parameter of distribution is given by $\delta=c_{2}%
\cos^{-2}(u+c_{1}),$ where $c_{1},c_{2}=const.$ The curve $\Gamma$ is planar.
In particular in the second case, the affine normal image of $\Phi$ is an
improper sphere and the curve $\Gamma$ is a straight line.
\end{proposition}

%-----------------------------------------------------------------------------------------------------------

We consider a conoidal surface $\Phi$ and a support function $q$ of the form
(\ref{2.7a}). For the curvature of the planar curve $\Gamma$ we find

\begin{equation}
k=\frac{Aff\
%TCIMACRO{\U{b4}}%
%BeginExpansion
\acute{}%
%EndExpansion
\
%TCIMACRO{\U{b4}}%
%BeginExpansion
\acute{}%
%EndExpansion
-2\delta^{2}Af\
%TCIMACRO{\U{b4}}%
%BeginExpansion
\acute{}%
%EndExpansion
\ ^{2}+Bff\
%TCIMACRO{\U{b4}}%
%BeginExpansion
\acute{}%
%EndExpansion
+Cf^{2}}{D}, \label{2.24}%
\end{equation}
where%
\[
A=\delta^{2}\left(  2\delta \delta%
%TCIMACRO{\U{b4}}%
%BeginExpansion
\acute{}%
%EndExpansion
\
%TCIMACRO{\U{b4}}%
%BeginExpansion
\acute{}%
%EndExpansion
-3\delta \
%TCIMACRO{\U{b4}}%
%BeginExpansion
\acute{}%
%EndExpansion
\ ^{2}-4\delta^{2}\right)  ,
\]%
\[
B=2\delta \left(  6\delta \delta \,%
%TCIMACRO{\U{b4}}%
%BeginExpansion
\acute{}%
%EndExpansion
\  \delta \
%TCIMACRO{\U{b4}}%
%BeginExpansion
\acute{}%
%EndExpansion
\
%TCIMACRO{\U{b4}}%
%BeginExpansion
\acute{}%
%EndExpansion
-4\delta^{2}\delta \
%TCIMACRO{\U{b4}}%
%BeginExpansion
\acute{}%
%EndExpansion
-6\delta \
%TCIMACRO{\U{b4}}%
%BeginExpansion
\acute{}%
%EndExpansion
\ ^{3}-\delta^{2}\delta \,%
%TCIMACRO{\U{b4}}%
%BeginExpansion
\acute{}%
%EndExpansion
\,%
%TCIMACRO{\U{b4}}%
%BeginExpansion
\acute{}%
%EndExpansion
\,%
%TCIMACRO{\U{b4}}%
%BeginExpansion
\acute{}%
%EndExpansion
\, \right)  ,
\]%
\[
C=4\delta^{4}+7\delta^{2}\delta \,%
%TCIMACRO{\U{b4}}%
%BeginExpansion
\acute{}%
%EndExpansion
\,^{2}+6\delta \,%
%TCIMACRO{\U{b4}}%
%BeginExpansion
\acute{}%
%EndExpansion
\,^{4}-2\delta^{3}\delta \,%
%TCIMACRO{\U{b4}}%
%BeginExpansion
\acute{}%
%EndExpansion
\,%
%TCIMACRO{\U{b4}}%
%BeginExpansion
\acute{}%
%EndExpansion
-6\delta \delta \,%
%TCIMACRO{\U{b4}}%
%BeginExpansion
\acute{}%
%EndExpansion
\,^{2}\delta \,%
%TCIMACRO{\U{b4}}%
%BeginExpansion
\acute{}%
%EndExpansion
\,%
%TCIMACRO{\U{b4}}%
%BeginExpansion
\acute{}%
%EndExpansion
+\delta^{2}\delta \,%
%TCIMACRO{\U{b4}}%
%BeginExpansion
\acute{}%
%EndExpansion
\, \delta \,%
%TCIMACRO{\U{b4}}%
%BeginExpansion
\acute{}%
%EndExpansion
\,%
%TCIMACRO{\U{b4}}%
%BeginExpansion
\acute{}%
%EndExpansion
\,%
%TCIMACRO{\U{b4}}%
%BeginExpansion
\acute{}%
%EndExpansion
,
\]%
\[
D=\left[  \delta^{2}\left(  \delta \,%
%TCIMACRO{\U{b4}}%
%BeginExpansion
\acute{}%
%EndExpansion
\,f-2\delta f\,%
%TCIMACRO{\U{b4}}%
%BeginExpansion
\acute{}%
%EndExpansion
\, \right)  ^{2}+\left[  2\left(  \delta^{2}+\delta \,%
%TCIMACRO{\U{b4}}%
%BeginExpansion
\acute{}%
%EndExpansion
\,^{2}\right)  f-\delta \left(  \delta \,%
%TCIMACRO{\U{b4}}%
%BeginExpansion
\acute{}%
%EndExpansion
\,f\,%
%TCIMACRO{\U{b4}}%
%BeginExpansion
\acute{}%
%EndExpansion
\,+\delta \,%
%TCIMACRO{\U{b4}}%
%BeginExpansion
\acute{}%
%EndExpansion
\,%
%TCIMACRO{\U{b4}}%
%BeginExpansion
\acute{}%
%EndExpansion
\,f\right)  \right]  ^{2}\right]  ^{3/2}.
\]
The curve $\Gamma$ is a straight line if and only if the function $f(u)$
fulfills the differential equation%
\[
Aff\
%TCIMACRO{\U{b4}}%
%BeginExpansion
\acute{}%
%EndExpansion
\,%
%TCIMACRO{\U{b4}}%
%BeginExpansion
\acute{}%
%EndExpansion
-2\delta^{2}Af\,%
%TCIMACRO{\U{b4}}%
%BeginExpansion
\acute{}%
%EndExpansion
\,^{2}+Bff\,%
%TCIMACRO{\U{b4}}%
%BeginExpansion
\acute{}%
%EndExpansion
+Cf^{2}=0.
\]
The existence of ruled surfaces, whose Laplace normal image indeed degenerates
into a straight line, whenever the support function is of the form
(\ref{2.7a})$,$ is guaranteed by the following examples.\medskip

%-----------------------------------------------------------------------------------------------------------

\noindent 1. Let $\Phi$ be a conoidal surface with
\[
\delta=const.\neq0.
\]
From
(\ref{2.24}) we obtain for the curvature of $\Gamma$%
\[
k=\frac{\left \vert \delta \left(  -ff\
%TCIMACRO{\U{b4}}%
%BeginExpansion
\acute{}%
%EndExpansion
\
%TCIMACRO{\U{b4}}%
%BeginExpansion
\acute{}%
%EndExpansion
+2f\
%TCIMACRO{\U{b4}}%
%BeginExpansion
\acute{}%
%EndExpansion
\ ^{2}+f^{2}\right)  \right \vert }{\left(  f\
%TCIMACRO{\U{b4}}%
%BeginExpansion
\acute{}%
%EndExpansion
\ ^{2}+f^{2}\right)  ^{3/2}}.
\]
Therefore $\Gamma$ is a straight line if and only if%
\[
ff\
%TCIMACRO{\U{b4}}%
%BeginExpansion
\acute{}%
%EndExpansion
\
%TCIMACRO{\U{b4}}%
%BeginExpansion
\acute{}%
%EndExpansion
-2f\
%TCIMACRO{\U{b4}}%
%BeginExpansion
\acute{}%
%EndExpansion
\ ^{2}-f^{2}=0,
\]
whose general solution, after a reparametrization, is%
\[
f(u)=\frac{c}{\cos u},\quad c=const.\neq0.
\]
Thus, \textit{the Laplace normal image of a conoidal surface with constant
parameter of distribution degenerates into a straight line if and only if the
support function is given by}%
\[
q(u,v)=\frac{c}{\cos u\,w},\quad c=const.\neq0.\medskip
\]

%-----------------------------------------------------------------------------------------------------------

\noindent 2. Let $\Phi$ be a conoidal surface with
\[
\delta=\sin^{2}u.\]
From (\ref{2.24})
we obtain for the curvature of $\Gamma$%
\[
k=4\sqrt{2}\sin^{4}u\frac{\left \vert \left(  6-5\sin^{2}u\right)  f^{2}%
+2\sin^{2}u\ f\
%TCIMACRO{\U{b4}}%
%BeginExpansion
\acute{}%
%EndExpansion
\ ^{2}-f\left(  3\sin2u\ f\
%TCIMACRO{\U{b4}}%
%BeginExpansion
\acute{}%
%EndExpansion
+\sin^{2}u\ f\
%TCIMACRO{\U{b4}}%
%BeginExpansion
\acute{}%
%EndExpansion
\
%TCIMACRO{\U{b4}}%
%BeginExpansion
\acute{}%
%EndExpansion
\  \right)  \right \vert }{\left \vert \left(  5\cos2u+13\right)  f^{2}%
-6\sin2u\ ff\
%TCIMACRO{\U{b4}}%
%BeginExpansion
\acute{}%
%EndExpansion
+2\sin^{2}u\ f\
%TCIMACRO{\U{b4}}%
%BeginExpansion
\acute{}%
%EndExpansion
\ ^{2}\right \vert ^{3/2}}.
\]
The curve $\Gamma$ is a straight line if and only if%
\[
f\left(  3\sin2u\ f\
%TCIMACRO{\U{b4}}%
%BeginExpansion
\acute{}%
%EndExpansion
+\sin^{2}u\ f\
%TCIMACRO{\U{b4}}%
%BeginExpansion
\acute{}%
%EndExpansion
\
%TCIMACRO{\U{b4}}%
%BeginExpansion
\acute{}%
%EndExpansion
\  \right)  +\left(  5\sin^{2}u-6\right)  f^{2}-2\sin^{2}u\ f\
%TCIMACRO{\U{b4}}%
%BeginExpansion
\acute{}%
%EndExpansion
\ ^{2}=0,
\]
whose general solution is%
\[
f(u)=\frac{c\sin^{3}u}{\cos2u},\quad c=const.,\neq0.
\]
Therefore, \textit{the Laplace normal image of a conoidal surface with
parameter of distribution given by }$\delta=\sin^{2}u$ \textit{degenerates
into a straight line if and only if the support function is given by}%
\[
q(u,v)=\frac{c\sin^{3}u}{\cos2u\,w},\quad c=const.\neq0.\medskip
\]

%-----------------------------------------------------------------------------------------------------------

\section{Ruled surfaces, whose Laplace normal image is a ruled surface}

In this paragraph we deal with non conoidal ruled surfaces, whose relative
normals $\bar{y}$ along each ruling \emph{lie on the corresponding asymptotic plane}.
On account of (\ref{1.6b}) it is easy to see, that this occur if and only if%
\[
vq+w^{2}q_{/2}=0,
\]
which gives that the support function is of the form (\ref{2.7a}), i.e.%

\[
q=\frac{f(u)}{w},
\]
where $f(u)$ is an arbitrary nonvanishing $C^{2}$-function.\medskip

\noindent Thus,
\[
\Phi^{\ast}:\bar{L}=\frac{f}{\delta}\bar{n}+\frac{f\left(  2\kappa v+\delta \,%
%TCIMACRO{\U{b4}}%
%BeginExpansion
\acute{}%
%EndExpansion
\, \right)  }{2\delta^{2}}\bar{e}
\]
is a parametrization of the Laplace normal image $\Phi^{\ast}$ of $\Phi$. We
see immediately, that $\Phi^{\ast}$ is a ruled surface too, whose rulings are
parallel to those of $\Phi$. So, $\Phi$ and $\Phi^{\ast}$ possess common
conical curvature and common Kruppa moving frame $\left \{  \bar{e},\bar
{n},\bar{z}\right \}  .$\medskip

We consider the directrix
\begin{equation}
\Gamma^{\ast}:\quad \bar{r}^{\ast}=\frac{f}{\delta}\bar{n} \label{2.51b}%
\end{equation}
of $\Phi^{\ast}$. Then, it is easy to prove the following

%-----------------------------------------------------------------------------------------------------------

\begin{proposition}
The tangents of the line of striction $\Sigma$ of $\Phi$ and the tangents of
the curve $\Gamma^{\ast}$ of $\Phi^{\ast}$ at corresponding points
are\smallskip \newline a) parallel if and only if
\begin{equation}
f=c \delta,\text{where }c=const.,
\end{equation}
and
\begin{equation}
\kappa \lambda+1=0,
\end{equation}
i.e., the support function is given by
\begin{equation}
q=\frac{c\, \delta}{w},
\end{equation}
and the line of striction $\Sigma$ of $\Phi$ is a line of curvature.\smallskip
\newline b) orthogonal if and only if
\begin{equation}
\kappa=\lambda,
\end{equation}
i.e., the line of striction $\Sigma$ of $\Phi$ is an asymptotic line.
\end{proposition}

A parametrization of the line of striction $\Sigma^{\ast}$ of $\Phi^{\ast}$
can be found on account of the relation $\bar{s}^{\ast}=\bar{r}^{\ast}%
-\langle \bar{r}^{\ast}\,%
%TCIMACRO{\U{b4}}%
%BeginExpansion
\acute{}%
%EndExpansion
,\bar{e}\,%
%TCIMACRO{\U{b4}}%
%BeginExpansion
\acute{}%
%EndExpansion
\, \rangle \, \bar{e}.$ We obtain%
\begin{equation}
\bar{s}^{\ast}=\frac{f}{\delta}\bar{n}-\left(  \frac{f}{\delta}\right)
%TCIMACRO{\U{b4}}%
%BeginExpansion
\acute{}%
%EndExpansion
\  \, \bar{e}. \label{2.52}%
\end{equation}
From (\ref{2.51b}) and (\ref{2.52}) we have:\medskip

%-----------------------------------------------------------------------------------------------------------

\begin{proposition}
\textit{The curve }$\Gamma^{\ast}$\textit{ coincides with the line of striction
}$\Sigma^{\ast}$\textit{ of }$\Phi^{\ast}$\textit{ if and only if the support
function is given by }%
\[
q=\frac{c\, \delta}{w},\mathit{\ }\text{where}\mathit{\ }c=const\mathit{.}%
\]

\end{proposition}

From (\ref{2.52}) we find%

\[
\bar{s}^{\ast}\acute{}=-\left[  \left(  \frac{f}{\delta}\right)
%TCIMACRO{\QATOP{\acute{}~\acute{}}{{}} }%
%BeginExpansion
\genfrac{}{}{0pt}{}{\acute{}~\acute{}}{{}}
%EndExpansion
+\frac{f}{\delta}\right]  \  \bar{e}+\frac{\kappa \,f}{\delta}\bar{z}.
\]
Consequently (see (\ref{1.23})), the fundamental invariants of $\Phi^{\ast}$
are the following%

\begin{equation}
\delta^{\ast}=\frac{\kappa \,f}{\delta},\quad \kappa^{\ast}=\kappa,\quad
\lambda^{\ast}=\frac{-\delta}{\kappa f}\left[  \left(  \frac{f}{\delta
}\right)
%TCIMACRO{\QATOP{\acute{}~\acute{}}{{}}}%
%BeginExpansion
\genfrac{}{}{0pt}{}{\acute{}~\acute{}}{{}}%
%EndExpansion
+\frac{f}{\delta}\right]  . \label{2.54}%
\end{equation}
From (\ref{2.54}) we obtain a series of results, which are contained in the following
%-----------------------------------------------------------------------------------------------------------

\begin{proposition}
(a) The Euclidean normals of $\Phi$ and $\Phi^{\ast}$ are parallel if and
only if%
\[
q=c~q_{_{AFF}},\text{ where }c=const.
\]
(b) The Laplace normal image $\Phi^{\ast}$ of $\Phi$ is orthoid $(\lambda^{\ast}=0)$ if and only if
the support function is given by
\[
q=\frac{\delta \left(  c_{1}\cos u+c_{2}\sin u\right)  }{w},\text{ where }%
c_{1},c_{2}=const.
\]
(c) The line of striction $\Sigma^{\ast}$ of $\Phi^{\ast}$ is an asymptotic
line $\left(  \kappa^{\ast}=\lambda^{\ast}\right)  $ if and only if
\begin{equation}
\delta \left(  \frac{f}{\delta}\right)  \,
\genfrac{}{}{0pt}{}{\acute{}~\acute{}}{{}}
%EndExpansion
+f\, \left(  \kappa^{2}+1\right)  =0. \label{2.55}
\end{equation}
(d) The line of striction $\Sigma^{\ast}$ of $\Phi^{\ast}$ is a line of
curvature\ $\left(  1+\kappa^{\ast}\lambda^{\ast}=0\right)  $ if and only if
the support function is given by
\[
q=\frac{\delta \left(  c_{1}+c_{2}u\right)  }{w}\text{, where }c_{1}%
,c_{2}=const.
\]
(e) $\Phi$ and $\Phi^{\ast}$ are congruent if and only if the support function
is given by
\[
q=\frac{\delta^{2}}{\kappa \,w}%
\]
and the fundamental invariants of $\Phi$ are associated with the relation
\[
\lambda=-\frac{1}{\delta}\left(  \frac{\delta}{\kappa}\right)  \,%
\genfrac{}{}{0pt}{}{\acute{}~\acute{}}{{}}%
%EndExpansion
-\frac{1}{\kappa}.
\]
(f) The Laplace normal image $\Phi^{\ast}$ is an Edlinger surface
$(\delta^{\ast}=const., \kappa^{\ast}\lambda^{\ast}+1=0)$ if and only if the
support function and the conical curvature are given by
\[
q=\frac{\delta \left(  c_{1}+c_{2}u\right)  }{w},\quad \kappa=\frac{c_{3}}%
{c_{1}+c_{2}u}\text{, where }c_{1},c_{2},c_{3}=const.
\]
respectively.
\end{proposition}
%-----------------------------------------------------------------------------------------------------------

The following example guarantees, that there exist ruled surfaces, which realize case (c) of the above proposition, i.e. ruled surfaces, such that the line of striction $\Sigma^{\ast}$ of their Laplace normal image $\Phi^{\ast}$ is an asymptotic line: For
\[
\kappa=c_{1}=const.
\]
we obtain from (\ref{2.55})
\[
f=\delta \left( c_{3}\cos c_{2}u+c_{4}\sin c_{2}u\right)  ,\quad c_{3}%
,c_{4}=const.,
\]
where $c_{2}=\sqrt{c_{1}^{2}+1}$. Consequently, when the support function is of the form
\[
q=\frac{\delta \left( c_{3}\cos c_{2}u+c_{4}\sin c_{2}u\right)  }{w}   ,\quad c_{3}%
,c_{4}=const.,\quad c_{2}=\sqrt{c_{1}^{2}+1},
\]
then, the line of striction $\Sigma^{\ast}$ of $\Phi^{\ast}$ is an asymptotic
line.
%-------------------------------------------------------------------------------------------------

\end{document}